\documentclass{article}
\usepackage{amsmath}
\usepackage{amssymb}
\usepackage{latexsym}
\usepackage{amsthm}
\usepackage{mathrsfs} 
\usepackage{wasysym}
\usepackage{fancyhdr}
\usepackage{amsfonts}
\usepackage{amssymb}
\usepackage{epsfig}
\usepackage{epstopdf}
\usepackage[normalem]{ulem}
\usepackage{eucal}

\usepackage{tkz-graph}
\usetikzlibrary{arrows,decorations.markings}
\usepackage{tkz-euclide}
\usetkzobj{all}

\newtheorem{thm}{Theorem}[section]

\newtheorem{lem}[thm]{Lemma}

\newtheorem{obs}[thm]{Observation}

\title{Connectivity for Kite-Linked Graphs}

\author{Chris Stephens\thanks{Department of Mathematical Sciences, Middle Tennessee State University, Murfreesboro, TN 37132. Email:\texttt{chris.stephens@mtsu.edu}.}\; and Dong Ye\thanks{Corresponding author. Department of Mathematical Sciences, Middle Tennessee State University, Murfreesboro, TN 37132. Email: \texttt{dong.ye@mtsu.edu}. Partially supported by a grant from Simons Foundation (no 359516).}}

\date{October 19 2019}

\begin{document}
\maketitle
\begin{abstract}
For a given graph $H$, a graph $G$ is {\em $H$-linked} if, for
every injection $\varphi: V(H) \to V(G)$, the graph $G$ contains a subdivision of $H$
with $\varphi(v)$ corresponding to $v$, for each $v\in V(H)$. Let $f(H)$ be the minimum integer $k$ such that every 
$k$-connected graph is $H$-linked. Among graphs $H$ with at least four vertices, the exact value $f(H)$ is only know when $H$ is a path with four vertices or a cycle with four vertices.
A \emph{kite} is graph obtained from $K_4$ by deleting two adjacent edges, i.e.,  a triangle together with a pendant edge. Recently, Liu, Rolek and Yu proved that every $8$-connected graph is kite-linked. The exact value of $f(H)$ when $H$ is the kite remains open. In this paper, we settle this problem by showing that every 7-connected graph is kite-linked.  
\medskip


\end{abstract}

\section{Introduction}
All graphs in this paper are simple and finite. A graph $G$ is $k$-connected if $G$ remains connected after deleting at most $k-1$ vertices of $G$. The well-known Menger's Theorem states that a graph $G$ is $k$-connected if and only if there exists $k$ disjoint paths between any pair of
disjoint vertex subsets of size at least $k$. A graph $G$ is {\em $k$-linked} if, for any given $2k$ distinct vertices $s_1, t_1, s_2, t_2, \ldots, s_k,t_k$, there are $k$ disjoint paths $P_1,\ldots, P_k$ such that each $P_i$ joins $s_i$ and $t_i$. Clearly, a $k$-linked graph is $k$-connected, but not vice versa.  
 Let $f(k)$ be the minimum 
value $t$ such that every $t$-connected graph is $k$-linked. The existence of the function $f(k)$ follows from a series of papers
by Larman and Mani \cite{LM}, Jung \cite{Jun} and Mader \cite{Mad}. The first linear bound $f(k)\le 22k$ was obtained by Bollob\'as and Thomason \cite{BT},  which was improved to $f(k)\le 12k$  by Kawarabayashi, Kostochka and G. Yu \cite{KKY} and to  $ f(k)\le 10k$ by Thomas and Wollan \cite{TW05}.

For 2-linked graphs,  Jung \cite{Jun} showed that a 4-connected non-planar graph is 2-linked and established that $f(2)=6$. A characterization for non-2-linked graphs is obtained independently by Seymour \cite{PS} and Thomassen \cite{CT}.  

\begin{thm}[Jung, \cite{Jun}]\label{thm:2-link}
Every 6-connected graph is 2-linked.
\end{thm}

The value $f(k)$ remains unknown for all $k\ge 3$. For 3-linked graphs, Thomas and Wollan \cite{TW08} proved that every 6-connected graph $G$ with at least $5|V(G)|-14$ edges are 3-linked, which implies that every 10-connected graph is 3-linked. 

 Let $H$ be a given graph. A graph $G$ is {\em $H$-linked} if, for every injection $\varphi:V(H)\to V(G)$, the graph $G$ contains a subdivision of $H$ with $\varphi(v)$ corresponding to $v$ for each $v\in V(H)$. The definition of $H$-linkage generalizes $k$-linkage. A graph is $k$-linked if $H$ is a union of $k$ independent edges, i.e., $k$ copies of $K_2$. For a given graph $H$, let $f (H)$ be the minimum 
 integer $\alpha$ such that every $\alpha$-connected graph is $H$-linked. 
 It follows immediately that  $f(H_1)\le f(H_2)$ if $H_1$ is a subgraph of $H_2$.
It is straightforward to see that $f(H)\ge |E(H)|$ and, in fact, this trivial lower bound is tight. For example, let $H$ be a star (i.e., $K_{1,t}$), 
then $f(H)=t$ which follows from Menger's Theorem. 
By the result of Thomas and Wollan~\cite{TW05} on the bound of $f(k)$, it follows that $f(H) \le 10 |E(H)|$. So, we have the following observation.

\begin{obs}
Let $H$ be a graph. Then $|E(H)|\le f(H)\le 10 |E(H)|$. 
\end{obs}

For small graphs other than $K_{1,t}$, $f(K_3)=3$ which follows from the classic result of Dirac \cite{Dirac} that every $k$-connected graph has a cycle through any $k$ vertices. It becomes nontrivial to determine the exact value $f(H)$ when $H$ has more
than three vertices. Ellingham, Plummer and G. Yu \cite{EPY} determined $f(P_4)=7$ where $P_4$ is a path with four vertices. So $f(H)\ge 7$ if $H\ne K_{1,3}$ has at least four vertices. In \cite{MWY}, McCarty, Wang and X. Yu show $f(C_4)= 7$, which was originally conjectured by Faudree \cite{Fau}.  
Very recently,  Liu, Rolek and G. Yu \cite{LRY} show that every 8-connected graph is kite-linked, where kite is a graph obtained from $K_4$ by deleting two adjacent edges. So $7\le f(H)\le 8$ if $H$ is the kite.  The proofs of these results where $H$ is $P_4$, $C_4$, or the kite in \cite{EPY, LRY, MWY} are based on fine structures of 3-planar graphs and obstructions developed in \cite{PS, XY03-1, XY03-2, XY03-3}.

The $H$-linkage problem has also been studied for special families of graphs. Goddard \cite{Go} showed that a 4-connected plane triangulation is $C_4$-linked. This theorem was extended to all surface triangulations by Mukae and Ozeki \cite{MO}. 
Ellingham, Plummer and G. Yu \cite{EPY} strengthen Goddard's result by proving that a 4-connected plane triangulation is $K_4^-$-linked, where $K_4^-$ denotes the graph obtained from $K_4$ by deleting one edge. 


The exact value $f(H)$ when $H$ is the kite remains open. In this paper, we settle this problem by proving the following result.

\begin{thm} \label{thm:main}
Every 7-connected graph is kite-linked. 
\end{thm}
 
It follows immediately from Theorem~\ref{thm:main}  that $f(H)=7$ if $H$ is the kite.

\section{Proof of Theorem~\ref{thm:main}}
In this section, we are going to prove Theorem~\ref{thm:main}. Before proceeding to the proof, we need some results from \cite{GPYZ, Per} and  \cite{LRY}.

Let $G$ be a graph and $S\subset V(G)$. Suppose $x\in V(G)- S$. A {\it $k$-fan between $x$ and $S$} consists of  $k$ internally disjoint paths joining $x$ and $k$ distinct vertices of $S$ such that each of the $k$ paths meets $S$ in exactly one vertex. 
An important result on connectivity due to Perfect~\cite{Per} states that, if a graph $G$ has a $k$-fan between $x$ and $S$, then any $(k-1)$-fan between $x$ and $S$ can be extended to a $k$-fan such that the end vertices of the $(k-1)$-fan are maintained in the $k$-fan. 
The following result is a special case of a generalized version of Perfect's Theorem proved in \cite{GPYZ} (Theorem 2.3).

\begin{thm} \label{thm:PT}
Let $G$ be a graph, and let $x\in V(G)$ and $S\subseteq V(G)\backslash \{x\}$ such that $G$ has a $k$-fan between $x$ and $S$.
If $S$ has a subset $T$ of size $k-t$ such that $G$ has a $(k-t)$-fan between $x$ and $T$,
then there exists a vertex subset $T'\subseteq S-T$ of size $t$ such that  $G$ has a $k$-fan between $x$ and $T\cup T'$.
\end{thm}

The following is a useful structure defined by Liu, Rolek and G. Yu in \cite{LRY} to find a kite-subdivision.

\medskip
\noindent {\bf Definition.} Let $x_1, x_2, x_3$ and $x_4$ be four give vertices of a graph $G$. A {\em $(x_1, x_2, x_3, x_4)$-flower} consisting of three cycles $C_1, C_2, C_3$ and three vertex-disjoint paths $P_1, P_2, P_3$ such that:

(1) $x_1\in C_1$, $x_3\in C_2$, $C_1\cap C_2=x_2$, and $C_3$ containing $x_4$ is disjoint from both $C_1$ and $C_2$;

(2) each path $P_i$, internally disjoint from the three cycles, joins the vertex $x_i$ to a vertex $v_i\in C_3$ for $i\in [3]$ such that $v_1, v_2, v_3, x_4$ are four different vertices appearing on $C_3$ in order. 

 \medskip

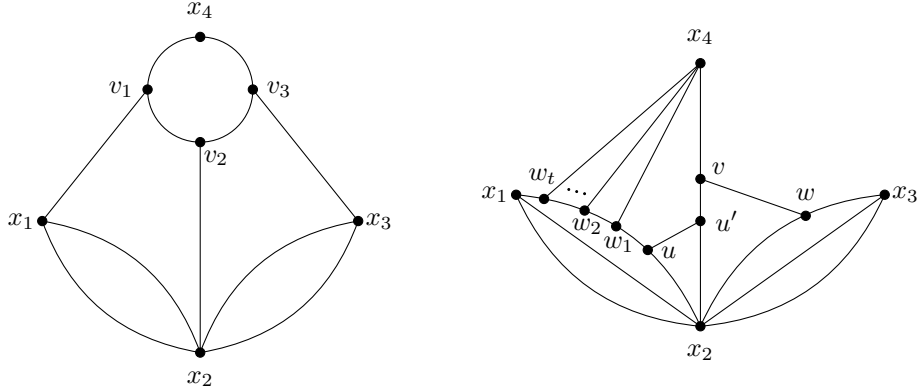
\begin{figure}[!hbtp]
\begin{center}
\begin{tikzpicture}[scale=.7]
 \tikzset{
  big arrow/.style={
    decoration={markings,mark=at position 1 with {\arrow[scale=1.3,#1]{>}}},
    postaction={decorate},
    shorten >=0.4pt},
  big arrow/.default=black}

\draw[] (0,3) circle (1); \filldraw[] (0,4) circle(2.5pt);\filldraw[] (0,2) circle(2.5pt);\filldraw[] (-1,3) circle(2.5pt);\filldraw[] (1,3) circle(2.5pt);
\filldraw[] (0, -2) circle(2.5pt); \filldraw[] (-3, .5) circle(2.5pt);  \filldraw[] (3, .5) circle(2.5pt);

 \draw[] (0,-2)  to [bend left=30]  (-3,.5) to [bend left=30] (0,-2);
 \draw[] (0, -2) to [bend left=30] (3, 0.5) to [bend left=30] (0,-2);
\draw[] (0,-2) -- (0, 2); \draw[] (3, .5)--(1, 3); \draw[] (-3,.5)--(-1, 3);

\node [] at (0, 4.5) {$x_4$}; \node [] at (1.5, 3) {$v_3$};  \node [] at (-1.5, 3) {$v_1$}; \node [] at (.3, 1.7) {$v_2$};
\node [] at (-3.4, .5) {$x_1$}; \node [] at (3.4, .5) {$x_3$};  \node [] at (0, -2.5) {$x_2$}; 

\begin{scope}[shift={(9.5,0.5)}]
\filldraw[] (0, 3) circle(2.5pt); \filldraw[] (0, -2) circle(2.5pt); \filldraw[] (-3.5, .5) circle(2.5pt);  \filldraw[] (3.5, .5) circle(2.5pt); 

\filldraw[] (-2.97, .42) circle(2.5pt); \filldraw[] (-1.6, -.1) circle(2.5pt); \filldraw[] (-2.2, .2) circle(2.5pt); \filldraw[] (0, 0) circle(2.5pt); \filldraw[] (-1.0, -.55) circle(2.5pt); \filldraw[] (2, .1) circle(2.5pt); \filldraw[] (0, .8) circle (2.5pt);

\filldraw[] (-2.2, .52)circle (.5 pt);  \filldraw[] (-2.35, .56)circle (.5 pt); \filldraw[] (-2.5, .6)circle (.5 pt);

 \draw[] (0,-2)  to [bend left=30]  (-3.5,.5) to [bend left=30] (0,-2)--(-3.5,.5);
 \draw[] (0, -2) to [bend left=30] (3.5, 0.5) to [bend left=30] (0,-2) --(3.5,.5);
\draw[] (0,-2) -- (0, 3); \draw[] (2.0, 0.1)--(0, .8); \draw[] (-3,.4)--(0, 3)--(-2.2, .2);\draw[] (0, 3)--(-1.6, -.1); \draw[] (0, 0)--(-1.0, -.55);

\node [] at (0, 3.5) {$x_4$}; \node [] at (-1.55, -.45) {$w_1$};  \node [] at (-2.15, -.1) {$w_2$};  \node [] at (-3.0, .85) {$w_t$};  \node [] at (-.6, -.6) {$u$}; \node [] at (.5, 0) {$u'$};
\node [] at (-3.9, .5) {$x_1$}; \node [] at (3.9, .5) {$x_3$};  \node [] at (0, -2.5) {$x_2$};  \node [] at (.35, 1) {$v$};  \node [] at (2, .5) {$w$}; 
\end{scope}

 \end{tikzpicture}
\caption{\small A flower (left) and an illustration for proof of Theorem~\ref{thm:main} (right).}\label{fig}
\end{center}
\end{figure}

\begin{lem}[Liu, Rolek and Yu, \cite{LRY}]\label{lem:flower}
If $G$ is a 7-connected graph containing a $(x_1, x_2, x_3, x_4)$-flower for some vertices $x_1, x_2, x_3, x_4\in V(G)$, then $G$ contains a kite-subdivision $K$ rooted at $x_1, x_2, x_3, x_4$ such that $x_4$ has degree 1 in $K$ and $x_2$ has degree 3 in $K$.  
\end{lem}

 Now, we are ready to prove Theorem~\ref{thm:main}. \medskip

\noindent{\bf Proof of Theorem~\ref{thm:main}.} Let $G$ be a 7-connected graph. Suppose to the contrary that $G$ is not kite-linked. 
So $G$ has four vertices $x_1, x_2, x_3$ and $x_4$ such that $G$ has no a kite-subdivision $K$ in which $x_4$ has degree 1 and $x_2$ has degree 3. 

By Menger's Theorem, $G$ has seven internally disjoint paths from $x_2$ to $\{x_1, x_3, x_4\}$ such that three of them, denoted by  $Q_i$ ($i\in [3]$), connect $x_2$ and $x_1$, another three, denoted by $R_i$ ($i\in [3]$), connect $x_2$ and $x_3$, and the remaining one $S$ connects $x_2$ and $x_4$. For convenience, let $Q=\cup_{i=1}^3 Q_i$ and $R=\cup_{i=1}^3 R_i$. 

Since $G$ is 7-connected, Menger's Theorem implies that there is a 7-fan between $x_4$ and $R\cup Q$. Note that, the path $S$ is internally disjoint from $R\cup Q$. Applying Theorem~\ref{thm:PT} to $x_4$ and $R\cup Q$, it follows that $G$ has a 7-fan between $x_4$ and $R\cup Q$,  which contains a path $P$ from $x_4$ to $x_2$. 
By symmetry,  assume that the 7-fan between $x_4$ and $R\cup Q$ contains at least three paths from $x_4$ to $Q$.  Let $W_1,...,W_t$ be all paths in the 7-fan from $x_4$ to $Q$, and let $w_1, ..., w_t$ be the landing vertices of these paths on $Q$. Then $t\ge 3$, and at least two of $w_1,..., w_t$ are different from $x_1$. Without loss of generality, assume that  $w_1, w_2\in V(Q)\backslash \{x_1, x_2\}$.  Let $W=\cup_{i=1}^t W_i$.

\medskip

\noindent {\bf Claim~1.} {\sl Every $(x_1, x_3)$-path intersects the path $P$.}

\medskip
\noindent{\em Proof of Claim~1.} Suppose $G$ has an $(x_1, x_3)$-path $P'$ which does not intersect $P$.  Then $P'$ has a subpath which joins a vertex $u$ from some $Q_i$ and a vertex $w$ from some $R_j$ and is internally disjoint from $R\cup Q$. Note that $u$ could be $x_1$ if $P'$ is internally disjoint from $Q$, and $w$ could be $x_3$ if $P'$ is internally disjoint from $R$. Without loss of generality, assume that $i=1$ if $u\ne x_1$ and $j=1$ if $w\ne x_3$. 
 Then $C=Q_2\cup R_2\cup R_1[x_3, w]\cup P'[w,u]\cup Q_1[u,x_1]$ together with the path $P$ is a desired kite-subdivision, a contradiction to the choice of $G$ as a counterexample. So Claim~1 follows. 
\medskip

Since $G$ is 7-connected, it follows from Theorem~\ref{thm:2-link} that $G$ is 2-linked. Hence, $G$ has two disjoint paths $L$ and $L'$ such that $L$ joins $x_1$ and $x_3$, and $L'$ joins $x_2$ and $x_4$. 
Note that $L$ does not contain $x_2$ and $x_4$. By Claim~1, $L$ intersects the path $P$. It follows that $P$ is not a single edge $x_2x_4$. 

Traverse along the path $L$ from $x_1$ to $x_3$. Let $u$ be the last intersecting vertex of 
$L$ and $Q\cup W$, and $v$ be the last intersecting vertex of $L$ and $P$, and $w$ be the first intersecting vertex of 
$L$ and $R$ after $v$. (Note that $u$ could be $x_1$ and $w$ could be $x_3$.) By Claim~1, $u, v$ and $w$ appear in order along $L$ in the direction from $x_1$ to $x_3$. 
By these assumptions, $L[v,w]$ is a path of  $L-E(Q\cup W\cup P\cup R)$, and the component  of $L-E(Q\cup W\cup P\cup R)$ containing $u$ contains a path $L[u,u']$ where $u'=L[u,u']\cap P$. 
Without loss of generality, we further assume that $w\in R_1$. Let $T=L[u,u']\cup P[u',v]\cup L[v,w]\cup R_1[w,x_3]\cup R_2$.

\medskip
\noindent {\bf Claim~2.} {\sl All vertices $w_1, ..., w_t$ belong to the same path from $Q_1, Q_2$ and $Q_3$. }
\medskip

\noindent{\em Proof of Claim~2.} Suppose to the contrary that at least two of $w_1,...,w_t$ belong to different paths. 
Without loss of generality, assume $w_1\in Q_1$ and $w_2\in Q_2$. 
 
If $u\in Q_i$ for some $i\in [3]$, at least one of $Q_1$ and $Q_2$ is different from $Q_i$, say $Q_1$.  Let $k\in [3]$ but $k\notin \{1, i\}$. Then the cycle $C=Q_k\cup Q_i[x_1,u]\cup T$ together with the path $W_1[x_4, w_1]\cup Q_1[w_1, x_2]$ form  a desired kite-subdivision, a contradiction. 

If $u\in W_i$ for some $i\in [t]$, then assume the landing vertex $w_i$ of $W_i$ on $Q$ belongs to $Q_j$ for some $j\in [3]$. Then
 at least one of $Q_1$ and $Q_2$ is different from $Q_j$, say $Q_1$. Let $k\in [3]$ but $k\notin \{i,j\}$. Then the cycle $C=Q_k\cup Q_j[x_1, w_i]\cup W_i[w_i,u] \cup T$ together with the path $W_1[x_4,w_1]\cup Q_1[w_1,x_2]$ form a desired kite-subdivision, a contradiction. This completes the proof of Claim~2. 
 \medskip

By Claim~2, assume all $w_1,...,w_t$ belong to $Q_1$. If $u$ belongs to $Q_2\cup Q_3$, then
$G$ has a desired kite-subdivision consisting of the cycle $C=Q_k\cup Q_j[x_1,u]\cup T$ and the path $W_1[x_4,w_1]\cup Q_1[w_1,x_2]$ where $j\in \{2,3\}$ and $k\notin\{1,j\}$, a contradiction. Hence $u$ belongs to either $Q_1$ or $W_k$ for some $k\in [t]$.

\medskip
\noindent {\bf Claim~3.} {\sl The vertex $w_i\in Q_1[x_1, u]$ for any $i\in [k]$  if $u\in Q_1$, and $w_i\in Q_1[x_1,w_k]$ for all $i\in [t]$ if $u\in W_k$.}
\medskip

\noindent{\em Proof of Claim~3.} If not, there exists a vertex, say $w_1$, failing the claim. 

If $u\in Q_1$, then $u\notin Q_1[w_1,x_2]$. Then $G$ has a desired kite-subdivision consisting of the cycle $C=Q_3\cup Q_1[x_1, u]\cup T$ and the path $W_1\cup Q_1[w_1,x_2]$, a contradiction. So assume that $u\in W_k$. Then $w_k\notin Q_1[ w_1, x_2]$. Then $G$ has a desired kite-subdivision consisting of the cycle $C=Q_3\cup Q_1[x_1, w_k]\cup W_k[w_k,u]\cup T$ and the path $W_1\cup Q_1[w_1, x_2]$, a contradiction again. 
This completes the proof of Claim~3. 
\medskip

By Claim~2, we may assume that $w_1, \ldots, w_t$ appear on $Q_1$ in order in the direction from $x_2$ to $x_1$. By Claim~3, it follows that $u, w_1, \ldots, w_t$ appear on $Q_1$ in order if $u\in Q_1$, and $u\in W_1$ if $u\notin Q_1$. (See Figure~\ref{fig}.)  

First, assume that $u\in Q_1$. Note that $u$ could be $w_1$, and $w_t$ could be $x_1$. Since $t\ge 3$ and $w_2$ is different from $w_1$ and $w_t$, it holds that $w_2\notin \{u, x_1\}$. 
Then let $C_3=W_2\cup Q_1[w_2, u]\cup L[u,u']\cup P[x_4, u']$, $C_1=Q_2\cup Q_3$, 
$C_2=R_2\cup R_3$, $P_1=Q_1[x_1, w_2]$, $P_2=Q_1[x_2, u]$, $P_3=R_1[x_3, w]\cup L[v,w]$ if $v\in P[x_4,u']$ 
and $P_3=R_1[x_3, w]\cup L[w,v]\cup P[v,u']$ otherwise. Then $x_4, w_2, u$ and $v$ if $v\in P[x_4, u']$ (or $u'$ if $v\in P(u', x_2)$) appear on $C_3$ in order. 
So $G$ has an $(x_1,x_2,x_3,x_4)$-flower. By Lemma~\ref{lem:flower}, $G$ has a 
desired kite-subdivision, a contradiction.

So, in the following, assume that $u\in W_1$.  Then $G$ has an $(x_1,x_2,x_3,x_4)$-flower: 
 $C_3=W_2\cup Q_1[w_2, w_1]\cup W_1[w_1,u]\cup L[u,u']\cup P[u',x_4]$, and $C_1, C_2, P_1,P_2,P_3$ are the same as described above. 
By Lemma~\ref{lem:flower}, $G$ has a desired kite-subdivision, a contradiction. This completes the proof. 
 \qed 

\end{document}